\documentclass[12pt]{article}
\usepackage{amsmath}
\usepackage{amssymb}
\usepackage{amscd}
\usepackage[dvips]{graphics}
\usepackage{amsfonts}

\newcommand{\refm [1]} {(\ref{#1})}

\newcounter{theorem}

\newtheorem{theorem}{Theorem}

\newcommand{\la}{\label}

\newcommand{\en}{\end{equation}}

\def\qed{\mbox{\rule{0.5em}{0.5em}}}

\def\BC{{\Bbb C}}

\def\qed{\mbox{\rule{0.5em}{0.5em}}}

\def\la{\label}

\def\del1{{\delta_n^{(1)}}}

\newcommand{\ignore}[1]{}

\def\BC{{\Bbb C}}

\newcommand{\be}{\begin{equation}}
\newcommand{\eu}{\end{equation}}
\newcommand{\ber}{\begin{eqnarray}}
\newcommand{\ena}{\end{eqnarray}}

\begin{document}
\title{ Explicit asymptotic formulae for sub exponentially growing multiplicative structures}
 \author{{\bf Boris L. Granovsky}
\thanks{E-mail: mar18aa@techunix.technion.ac.il} \\
Department of Mathematics, Technion-Israel Institute of Technology,\\
Haifa, 32000,Israel.}\maketitle


\begin{abstract}
 It is derived the explicit asymptotic expression in $n$ for the coefficient    $c_n$ of the generating function for  multiplicative structures with sub exponential rate of growth of $c_n,$ as $n\to\infty$.
\end{abstract}
\maketitle
{\bf Keywords:} Partitions, Generating function, Multiple poles, Asymptotic formula, Sub exponential growth, Local limit theorem.\\
{\bf MSC :}
05A17,11P82
\vskip 1cm

{\bf I. Mathematical background.}

Our starting point is the following  mathematical formalism adopted from \cite{Nes} and from the preceding  works \cite{multi},\cite{ext} with Dudley Stark.

Let the function $f(z),\ \vert z\vert\le 1,\ z\in \BC$ be analytic on the unit circle and  be of the following multiplicative form
\be f(z)=\prod_{j\ge 1} (S(a_jz^j))^{b_j}, \ \vert z\vert\le 1,\ z\in \BC, \la{bhx}\eu where all Taylor coefficients of the base function $S$ are $\ge 0,$  $S(0)=1$ and $$0<a_j\le 1, \ b_j\ge 0,\ j\ge 1$$ are two given sequences of  parameters.
This amounts to say that $f$ is the generating function for the sequence $c_n\ge 0,\ n\ge 0,\  c_0=1,$ in the sense that
$$f(z)=\sum_{n\ge 0} c_nz^n,\ \vert z\vert\le 1. $$
Moreover, due to the multiplicative form \refm[bhx] of $f$, the  Taylor coefficient $c_n$ of $f$ can be viewed as the number of the associated decomposable combinatorial structures of size $n$.
Varying the base function $S$ in \refm[bhx], a great scope of decomposable combinatorial structures can be encompassed. For example, the function  $S(z)=(1-z)^{-1}$ induces  integer partitions with each part $j$ weighted with $b_j$ and scaled with $a_j$ (more explanation on the meaning of the two parameters can be found in \cite{ext}).\\
It is appropriate to note that the above setting can be  spelled also in the language of arithmetic semigroups. For more details and references see the recent paper \cite{min}.

In \cite{Nes} there were formulated necessary and sufficient conditions providing sub exponential growth of $c_n,\ n\to \infty,$ given by the following asymptotic formula:
\be c_n\sim \frac{\delta_n^{\frac{\rho_r}{2}+1}}{\sqrt{2\pi}}\exp\left(\sum_{l=0}^r h_l\delta_n^{-\rho_l}
-A_0\log\delta_n + \Delta(\delta_n) + n\delta_n\right), \ n\to\infty.\la{expgr}\eu
We recall below the ingredients of the formula \refm[expgr].
\begin{itemize}
\item
$0<\rho_1<\ldots<\rho_r, r\ge 1$ are simple
poles of the function  $\Gamma(s)D(s), \ s\in \BC $, with respective residues $h_l>0,\ l=1,\ldots,r,$ where  \be\label{Dirdef}
D(s)=\sum_{k=1}^\infty \Lambda_k k^{-s},
\en
is the Dirichlet generating function for the Taylor coefficients of the function
\be
\log f(z) = \sum_{k\ge 1}b_k \log S(a_kz^k)=\sum_{k\ge 1} \Lambda_k z^k,\quad \vert z\vert<1.\la{lofl}\eu
The real coefficients $A_0,h_0$ are related to the pole  $\rho_0=0$\ of the function $\Gamma(s)D(s)$, which may be  either simple, if $D$ is analytic at $0$, or double, if $D$ has a simple pole at $0$,  namely:
$A_0= \lim_{s\to 0} (sD(s)),$ and, finally, $h_0= \Theta-\gamma A_0,\ \Theta=\lim_{s\to 0}( D(s)- A_0s^{-1})$.
\item
$\delta_n>0$ is the unique solution of the Khintchine equation
 \be \sum_{l=1}^r h_l\rho_l\delta_n^{-\rho_l-1}
 + A_0\delta_n^{-1} + D(-1) = n. \la{dltn}\eu

\item
$$\Delta(\delta_n)\to 0,\ n\to \infty$$ is the remainder term given by the expansion
\be \Delta(\tau)=\sum_{l\ge 1}\frac{(-1)^lD(-l)}{l!}\ \tau^l,\ \tau\in {\cal C}.\la{atrh,}\eu
\item
The out exponential term $\frac{\delta_n^{\frac{\rho_r}{2}\ + 1}}{\sqrt{2\pi}}$ comes from the normal local limit theorem which is a part of Khintchine's representation of $c_n$ (see the aforementioned references). It was proven in \cite{Nes} that for local limit theorem to hold it is sufficient that
the weights $\{b_j\}$
obey the condition
\be  \sum_{1\le k\le n,q\!\!\not| k} b_k\ge O(\log^2 n),\ \text{for\ any \ integer}\ q>1. \la{njkv}\eu

The conditions \refm[njkv] requires that for any $q>1$ the number of $b_j,\ j=1,\ldots,n$ divisible by $q$ should not be too large.

\end{itemize}
The purpose of the present paper is to derive  the explicit expression in $n$ of  the RHS of the asymptotic formula \refm[expgr]. In \cite{multi} the asymptotic formula for $c_n$ was expressed via powers of $n$ with coefficients defined recursively, i.e. non explicitly. The latter fact prompted the authors of \cite{plic} to seek for an explicit formula for the asymptotics of $c_n$.  They achieved the goal for a particular case of partitions with two specific poles $\rho_1=1,\rho_2=2$. We discuss their result in more details in the sequel of the present paper.

\vskip .5cm
{\bf II. Main theorem.}
\begin{theorem}
Let the two largest poles $\rho_{r-1},\rho_r$ obey the condition
\be 2\rho_{r-1}-\rho_r\le 0. \la{adrv}\eu
Then the following asymptotic formula holds:
$$c_n\sim \frac{1}{\sqrt{2\pi(\rho_r h_r)(\rho_r+1)}}(\rho_rh_r)^{\frac{\rho_r+2 -2A_0}{2(\rho_r +1)}}n^\kappa e^{Q}$$\be \exp{\Big((1+\rho_r)h_r(\rho_rh_r)^{-\frac{\rho_r}{\rho_r+1}}n^{\frac{\rho_r}{\rho_r+1}} +
\sum_{l=r-1}^1  h_l (\rho_rh_r)^{-\frac{\rho_l }{\rho_r +1}} n^{\frac{\rho_l }{\rho_r +1}}}\Big),\la{pjs}\eu
where we denoted
$$\kappa= (-\frac{\rho_r}{2}- 1 + A_0)(\rho_r +1)^{-1},$$
\be
Q=
\left\{
  \begin{array}{ll}
    h_0, & \hbox{if\ $2\rho_{r-1}-\rho_r<0$} \\
    h_0-(\rho_rh_r)^{-\frac{2\rho_{r-1}+1}{\rho_r+1}}
\frac{(\rho_{r-1}h_{r-1})^2}{2(\rho_r+1)}, & \hbox{if\ $2\rho_{r-1}-\rho_r=0$}
  \end{array}
\right. \la{trpk}\eu

\end{theorem}
{\bf Proof.}

{\bf Step 1: The solution of the Khintchine equation.}\\ Firstly, we denote $z_n=\delta_n^{-1},$  to rewrite the equation \refm[dltn] as \be \sum_{l=1}^r h_l\rho_lz_n^{\rho_l+1}
 + A_0z_n + D(-1) = n,\ z_n\to \infty,\ n\to \infty.\la{shftui}\eu
We will show that, under  the condition \refm[adrv], only the first two terms in the asymptotic expansion of $z_n,\ n\to \infty$ matter.\\
It follows from \refm[shftui] that $z_n\sim (\rho_rh_r)^{-\frac{1}{\rho_r +1}} n^{\frac{\rho_l}{\rho_r +1}},$ so that
\be z_n= (\rho_rh_r)^{-\frac{1}{\rho_r +1}} n^{\frac{\rho_r}{\rho_r +1}} + v_{1,n}\la{xckb}\eu
with $v_{1,n}=o(n^{\frac{\rho_r}{\rho_r +1}})$.
Then, substituting the above representation of $z_n$ into   \refm[shftui] and applying the binomial formula we get
$$ \rho_r h_r (\rho_r+1)(\rho_rh_r)^{-\frac{\rho_r}{\rho_r +1}} n^{\frac{\rho_r}{\rho_r +1}}v_{1,n}\sim
\rho_{r-1} h_{r-1}(\rho_rh_r)^{-\frac{\rho_{r-1}+ 1}{\rho_r +1}} n^{\frac{\rho_{r-1}+1}{\rho_r +1}},$$
which implies  $$v_{1,n}\sim M(\rho_rh_r)^{-\frac{\rho_{r-1}-\rho_r+1}{\rho_r +1}} n^{\frac{\rho_{r-1}-\rho_r+1}{\rho_r +1}}:= w_n,$$
where the constant $M=\frac{\rho_{r-1} h_{r-1}}{(\rho_r +1)\rho_r h_r}>0.$
In the same manner we have
$$z_n=(\rho_rh_r)^{-\frac{1}{\rho_r+1}} n^{\frac{1}{\rho_r +1}}+ w_n  +  v_{2,n},$$
with $v_{2,n}=o(v_{1,n}).$
Thus, $(\rho_rh_r)^{-\frac{\rho_l+1}{\rho_r+1}} n^{\frac{\rho_l+1}{\rho_r +1}}$ is the principal term in the binomial expansion, in $n\to \infty,$ of $z_n^{\rho_l+1},\ l=r,\ldots,1$. Next, we denote by
\be L_1= \sum_{l=r-1}^1 \rho_l h_l (\rho_rh_r)^{-\frac{\rho_l +1}{\rho_r +1}} n^{\frac{\rho_l +1}{\rho_r +1}} \la{srl}\eu  the sum of the above principal terms, except the principal term for $z_n^{\rho_r +1}$. We  also need the second and the third terms of the binomial  expansions of  $\rho_rh_rz_n^{\rho_r+1},$  denoted by $P_{1,r}$ and  $R_{1,r}$ respectively, and the same for the expansion of  $\rho_{r-1}h_{r-1}z_n^{\rho_{r-1}+1}$, denoted respectively by $P_{1,r-1}$ and $R_{1,r-1}$:
\be P_{1,r}= (\rho_r h_r)(\rho_r +1)(\rho_rh_r)^{-\frac{\rho_r}{\rho_r +1}} n^{\frac{\rho_r}{\rho_r +1}}(w_n + v_n^{(2)}),\la{,busv}\eu
\be R_{1,r}= (\rho_r h_r)\frac{(\rho_r +1)\rho_r}{2}(\rho_rh_r)^{-\frac{\rho_r -1}{\rho_r +1}} n^{\frac{\rho_r-1}{\rho_r +1}}(w_n + v_n^{(2)})^2=  \la{vpudv}\eu
$$(\rho_rh_r)^{-\frac{2\rho_{r -1}+2}{\rho_r +1}}\frac{\rho_r(\rho_{r-1}h_{r-1})^2}{2(\rho_r+1)}n^{\frac{2\rho_{r-1}-\rho_r +1}{\rho_r +1}}+ o(n^{\frac{2\rho_{r-1}-\rho_r +1}{\rho_r +1}}).$$
$$P_{1,r-1}= (\rho_{r-1} h_{r-1})(\rho_{r-1} +1)(\rho_{r}h_{r})^{-\frac{\rho_{r-1}}{\rho_r +1}} n^{\frac{\rho_{r-1}}{\rho_r +1}}(w_n + v_n^{(2)})=$$\be \frac{(\rho_{r-1}h_{r-1})^2(1+\rho_{r-1})}{1+\rho_r}(\rho_rh_r)^{-\frac{2\rho_{r-1} +2}{\rho_r +1}} n^{\frac{2\rho_{r-1}-\rho_r +1}{\rho_r +1}} + o(n^{\frac{2\rho_{r-1}-\rho_r +1}{\rho_r +1}}), \la{suet} \eu
\be R_{1,r-1}= (\rho_{r-1} h_{r-1})\frac{(\rho_{r-1} +1)\rho_{r-1}}{2}(\rho_rh_r)^{-\frac{\rho_{r-1} -1}{\rho_r +1}} n^{\frac{\rho_{r-1}-1}{\rho_r +1}}(w_n + v_n^{(2)})^2. \la{xdur} \eu
Finally, by $o_{1,r}$ and $o_{1,r-1}$ we denote the remainders of the two aforementioned expansions after the terms $R_{1,r}$ and $R_{1,r-1}$ respectively, while $O_{1,l}$ denotes the remainder of the expansion of $\rho_lh_lz_n^{\rho_l+1}, \ l=r-2,\ldots,1$ after the principal term.
By the preceding discussion,
$o_{1,r}=o(n^{\frac{2\rho_{r-1} - \rho_{r} +1}{\rho_r +1}}), o_{1,r-1}= o(n^{\frac{3\rho_{r-1} - 2 \rho_{r} +1}{\rho_r +1}})\ ,O_{1,l}=O(n^{\frac{\rho_l+ \rho_{r-1} - \rho_{r} +1}{\rho_r +1}}), \ l=r-2,\ldots,1$.
Invoking the above notations allows  to rewrite the equation \refm[shftui] in the following form:
$$L_1 + P_{1,r} + R_{1,r} + o_{1,r} + P_{1,r-1} + R_{1,r-1} + o_{1,r-1}  + \sum_{l=r-2}^1
O_{1,l} \ + $$\be A_0\Big((\rho_rh_r)^{-\frac{1}{\rho_r+1}} n^{\frac{1}{\rho_r +1}}+w_n + v_n^{(2)}\Big) + h_0 + D(-1)=0.\la{,phkv}\eu
{\bf Step 2: Asymptotic expansion of $\log c_n,\ n\to \infty.$}\\
Now our purpose will be to split, in a way similar to \refm[,phkv], the asymptotic expansion for $\log c_n,\ n\to \infty.$  Denoting by
$\bullet_2$ the quantities  corresponding to  $\bullet_1$ in  \refm[,phkv] we have from \refm[expgr],
$$ \log c_n= (-\frac{\rho_r}{2}- 1 + A_0)\log z_n - 1/2\log(2\pi)- 1/2\log(\rho_rh_r)-1/2 \log(\rho_r+1)   + $$$$
 \sum_{l=1}^r(1+\rho_l) h_lz^{\rho_l}_n + h_0+ A_0=   O(\log n) +$$
$$(1+\rho_r) h_r(\rho_rh_r)^{-\frac{\rho_r}{\rho_r+1}} n^{\frac{\rho_r}{\rho_r +1}} + L_2 +
P_{2,r}+ R_{2,r}+P_{2,r-1}+ o_{2,r} + R_{2,r-1} + o_{2,r-1}  + $$\be \sum_{l=r-2}^1 O_{2,l} + A_0 + h_0,
\la{na,sk} \eu
where for the sake of brevity we denoted $$O(\log n)=(\frac{-\rho_r}{2}- 1 + A_0)\log z_n - 1/2\log(\rho_rh_r)-1/2 \log(\rho_r+1)- 1/2\log(2\pi)=$$$$ (\frac{-\rho_r}{2}- 1 + A_0)\big(\frac{1}{\rho_r +1}\log(\rho_rh_r)+ \frac{1}{\rho_r +1} \log n\big) - 1/2\log(2\pi) - 1/2\log(\rho_rh_r)-1/2 \log(\rho_r+1) + \epsilon_n . $$ (Here the last equation follows from \refm[xckb]).

We now find the asymptotics of the terms $\bullet_2$ in the RHS of  \refm[na,sk].

$$L_2=\sum_{l=r-1}^1 (1+\rho_l)h_l(\rho_rh_r)^{-\frac{\rho_l}{\rho_r+1}} n^{\frac{\rho_l}{\rho_r +1}}=$$
$$\sum_{l=r-1}^1(1 + \rho_l^{-1})\rho_lh_l(\rho_rh_r)^{-\frac{\rho_l}{\rho_r+1}} n^{\frac{\rho_l}{\rho_r +1}}=$$
$$(\rho_rh_r)^{\frac{1}{\rho_r+1}}n^{-\frac{1}{\rho_r +1}}L_1 + \sum_{l=2}^1 h_l(\rho_rh_r)^{-\frac{\rho_l}{\rho_r+1}} n^{\frac{\rho_l}{\rho_r +1}},$$
where the last step follows from \refm[srl];
$$P_{2,r}= (1+\rho_r)\rho_r h_r(\rho_rh_r)^{-\frac{\rho_r -1}{\rho_r +1}} n^{\frac{\rho_r -1}{\rho_r +1}}(w_n + v_n^{(2)})=
$$$$(\rho_rh_r)^{\frac{1}{\rho_r +1}} n^{-\frac{1}{\rho_r +1}}P_{1,r},$$
by \refm[,busv];
$$R_{2,r}= (\rho_r+1) h_r \frac{\rho_r (\rho_r-1)}{2}(\rho_rh_r)^{-\frac{\rho_r -2}{\rho_r +1}} n^{\frac{\rho_r-2}{\rho_r +1}}(w_n + v_n^{(2)})^2=$$$$(1-\rho_r^{-1})(\rho_rh_r)^{\frac{1}{\rho_r +1}} n^{-\frac{1}{\rho_r+1}}R_{1,r},$$
by \refm[vpudv];
$$P_{2,r-1}= (1+\rho_{r-1})\rho_{r-1} h_{r-1}(\rho_rh_r)^{-\frac{\rho_{r-1} -1}{\rho_r +1}} n^{\frac{\rho_{r-1} -1}{\rho_r +1}}(w_n + v_n^{(2)})= (\rho_{r}h_{r})^{\frac{1}{\rho_r +1}} n^{-\frac{1}{\rho_r +1}}P_{1,r-1},$$
by \refm[suet];

$$R_{2,r-1}=(\rho_{r-1}+1) h_{r-1} \frac{\rho_{r-1} (\rho_{r-1}-1)}{2}(\rho_rh_r)^{-\frac{\rho_{r-1} -2}{\rho_r +1}} n^{\frac{\rho_{r-1}-2}{\rho_r +1}}(w_n + v_n^{(2)})^2=$$$$(1-\rho_{r-1}^{-1})(\rho_rh_r)^{\frac{1}{\rho_r +1}} n^{-\frac{1}{\rho_r+1}}R_{1,r-1},$$
by \refm[xdur];

$$o_{2,r}= O(n^{-\frac{1}{\rho_r+1}})o_{1,r}=o(n^{\frac{2\rho_{r-1}-\rho_r}{\rho_r+1}})\to 0,$$
$$o_{2,r-1}= O(n^{-\frac{1}{\rho_r+1}})o_{1,r-1}=o(n^{\frac{3\rho_{r-1} - 2 \rho_{r} }{\rho_r +1}})\to 0.$$
We also have, $$\rho_{r-1}^{-1}(\rho_rh_r)^{-\frac{1}{\rho_r +1}} n^{-\frac{1}{\rho_r+1}}R_{1,r-1}=$$
\be O(n^{-\frac{1}{\rho_r+1}}) O(n^{\frac{\rho_{r-1}-1}{\rho_r+1}})O(n^{\frac{2(\rho_{r-1}-\rho_r +1)}{\rho_r+1}})= O(n^{\frac{3\rho_{r-1}-2\rho_r}{\rho_r+1}})\to 0,\la{abhmss}\eu
by \refm[xdur],  the definition of $w_n$ and \refm[adrv].  Finally,
\be \rho_r^{-1}(\rho_rh_r)^{\frac{1}{\rho_r+1}}n^{-\frac{1}{\rho_r+1}}R_{1,r}= (\rho_rh_r)^{-\frac{2\rho_{r-1}+1}{\rho_r+1}}
\frac{(\rho_{r-1}h_{r-1})^2}{2(\rho_r+1)} n^{\frac{2\rho_{r-1}-\rho_r}{\rho_r+1}} +o(n^{\frac{2\rho_{r-1}-\rho_r}{\rho_r+1}}).\la{ktcru}\eu
Combining the above results, we derive
from \refm[expgr],$$\log c_n=$$ $$(1+\rho_r) h_r(\rho_rh_r)^{-\frac{\rho_r}{\rho_r+1}} n^{\frac{\rho_r}{\rho_r +1}} + L_2 +
P_{2,r}+ R_{2,r}+P_{2,r-1}+ o_{2,r} + R_{2,r-1} + o_{2,r-1}  + $$\be \sum_{l=r-2}^1 O_{2,l} + A_0 + h_0+ O(\log n),
 \eu
Applying now the above established  relationships between the quantities $\bullet_2$ and $\bullet_1,$ \refm[abhmss],\refm[ktcru] and the equation \refm[,phkv], we get
$$\log c_n= O(\log n) +
\sum_{l=1}^r (\rho_l+1)h_lz_n^{\rho_l}+ h_0 +  A_0 + \epsilon_n=$$$$ (-\frac{\rho_r}{2}- 1 + A_0)
\big(-\frac{1}{\rho_r +1}\log(\rho_rh_r)+ \frac{1}{\rho_r +1} \log n\big)+ $$$$(1+\rho_r) h_r(\rho_rh_r)^{-\frac{\rho_r}{\rho_r+1}} n^{\frac{\rho_r}{\rho_r +1}}+
\sum_{l=r-1}^1  h_l (\rho_rh_r)^{-\frac{\rho_l }{\rho_r +1}} n^{\frac{\rho_l }{\rho_r +1}} -A_0 +$$
$$-A_0(w_n + v_n^{(2)})(\rho_rh_r)^{\frac{1 }{\rho_r +1}}n^{-\frac{1 }{\rho_r +1}} + h_0 +  A_0 -$$
$$\rho_r^{-1}(\rho_rh_r)^{\frac{1}{\rho_r +1}} n^{-\frac{1}{\rho_r+1}}R_{1,r}
-\rho_{r-1}^{-1}(\rho_rh_r)^{\frac{1}{\rho_r +1}} n^{-\frac{1}{\rho_r+1}}R_{1,r-1} + O(n^{-\frac{1}{\rho_r+1}})
\Big( O_{1,l} + o_{1,r} + o_{1,r-1}\Big) + \epsilon_n=
$$
$$ (\frac{-\rho_r}{2}- 1 + A_0)
\big(-\frac{1}{\rho_r +1}\log(\rho_rh_r)+ \frac{1}{\rho_r +1} \log n\big)  - 1/2\log(2\pi)- 1/2\log(\rho_rh_r)-$$$$1/2 \log(\rho_r+1) + \epsilon_n + $$
$$(1+\rho_r) h_r(\rho_rh_r)^{-\frac{\rho_r}{\rho_r+1}} n^{\frac{\rho_r}{\rho_r +1}}+
\sum_{l=r-1}^1  h_l (\rho_rh_r)^{-\frac{\rho_l }{\rho_r +1}} n^{\frac{\rho_l }{\rho_r +1}} + h_0 - $$
$$\rho_r^{-1}
(\rho_rh_r)^{-\frac{2\rho_{r-1} +1}{\rho_r +1}}\frac{\rho_r(\rho_{r-1}h_{r-1})^2}{2(\rho_r+1)} n^{\frac{2\rho_{r-1}-\rho_r}{\rho_r+1}} +o(n^{\frac{2\rho_{r-1}-\rho_r}{\rho_r+1}})- O(n^{\frac{3\rho_{r-1}-2\rho_r}{\rho_r+1}}).$$
This leads to the stated asymptotic formula \refm[pjs]\qed

{\bf Two remarks.} (i) In view of the dichotomy in \refm[trpk], the relation $2\rho_{r-1}- \rho_r=0$
can be considered  as a critical one for the asymptotic behavior of $c_n,\ n\to \infty.$

(ii) The asymptotic formula \refm[pjs] derived under the condition \refm[adrv] includes only the first sum in the exponent of the general formula in Theorem 1 of \cite{multi}, since the second sum vanishes as $n\to\infty$, by  virtue of \refm[adrv].     In the case when the condition \refm[adrv] fails, the asymptotics of $c_n$ can be derived in the same way as in the present paper, by involving additional terms in the binomial expansions considered.

{\bf III. Corollaries}

We will demonstrate that the  formula \refm[pjs] encompasses  some known results.

 {\bf 1. Meinardus' theorem (1954).} This seminal formula was derived for weighted partitions with the Dirichlet generating function $D_b(s)$ for weights $b_j,\ j\ge 1$ having only one pole. We set  in \refm[pjs] $r=1$ and  recall  that for weighted partitions the function $D(s)$ is of the form $D(s)=\zeta(s+1)D_b(s),$
where $D_b(s)=\sum_{k\ge 1} \frac{b_k}{k^s}. $
The latter says that $h_0=(\Gamma(s+1)s\zeta(s+1)D_b(s))^\prime\vert_{s=0}= D_b^\prime(0),$ due to the fact that $\Gamma^\prime(1)=-\gamma$, while $h_l=A_l\zeta(\rho_l +1)\Gamma(\rho_l),\ l=1,\ldots,r,$ where $A_l>0,\ l=1,\ldots,r$ are poles of $D_b(s)$
and $A_0=D_b(0).$  Consequently, we have  $Q=h_0=D_b^\prime(0)$ in \refm[trpk], recovering the Meinardus' formula (see \cite{and},\cite{GSE}).\\
A particular case of the above model is the standard partitions for which $$\rho_r=1,\ D_b(s)=\zeta(s),\ D(s)=\zeta(s)\zeta(s+1),\ h_0=\zeta^\prime(0)= -\frac{1}{2} \log( 2\pi),$$$$D(0)=\zeta(0)=-\frac{1}{2}, \ h_r= \zeta(2)=\frac{\pi^2}{6}.$$
Substituting these data in \refm[pjs] gives the first term of the famous Hardy-Ramunajan expansion (1918) of $c_n$ (see e.g.\cite{and}).

{\bf 2. Partitions into roots. }

The subject of this example are  partitions  of an integer $n$ into integer parts of square roots of  the sequence of integers  $[\sqrt{1},
], [\sqrt{2}], \ldots , [\sqrt{l}],$ where  $l\ge 1$ is an integer.
Formally, the model is incorporated in the framework of our  Theorem 1, by setting
$$ D_b(s)=2\zeta(s-1)+\zeta(s),$$
because in the sequence $[\sqrt{1}],[\sqrt{2}],\ldots,[\sqrt{k}],\ldots$ a part $j$ is repeated $b_j=(j +1)^2-j^2= 2j+1$ times, so that
$$D(s)=\zeta(s+1)(2\zeta(s-1)+\zeta(s)).$$
Thus, $$r=2,\ \rho_2=2,\ \rho_{1}=1,$$ so that
 $$h_1= \zeta(2),h_2= 2\zeta(3), \ A_0= D_b(0)= -1/2 + 2\zeta(-1)= -1/2 -2 1/12=-2/3,$$ $$h_0=D_b^\prime(0)=\zeta^\prime(0)+ 2\zeta^\prime(-1)=-1/2\log(2\pi)+ 2 /12-2\log \gamma.$$
 Consequently,
$$Q= 1/6 - 1/2(4\zeta(3))^{-1}\zeta^2(2/6))=1/6- \frac{\zeta^2(2)}{24\zeta(3)}, $$
in accordance with \refm[trpk], since in the case considered $2\rho_{r-1}-\rho_r=0,$
while
$$ \kappa=(-2 - 2/3)1/3=- 8/9$$
and, finally, the out exponential term is equal to
$$\frac{1}{\sqrt{2\pi}\gamma^2}\frac{1}{\sqrt{2\pi 3(4\zeta(3))}}(4\zeta(3))^{(4+4/3)(6^{-1})}=
\frac{(4\zeta(3))^{7/18}}{\pi\gamma^2\sqrt{12}}.$$
The above data recover perfectly  the asymptotic formula in Theorem 2 in \cite{plic}, obtained by F.Luca and
D.Ralaivaosaona with the help of a combination of Meinardus' and the circle methods.  The circle method was
invented by Hardy and Ramanujan (1918) and since that time it has been used  by  generations of  researchers for asymptotic  enumeration. In the present paper, as well as in \cite{GSE}-\cite{ext} we used instead of the  circle method, the Khintchine approach, which we found more convenient for implementation.

{\bf 3. Partitions into congruent parts.}

For this model the generating function for the number of partitions is of the form
$$f(z)=\prod_{j\ge 1} \frac{1}{1-z^{aj +b}}, \ \vert z\vert<1,\ \text{where}\ a,b:gcd(a,b)=1, $$
which says that  $f$ counts  integer partitions into parts congruent to $b(mod\ a)$. To incorporate the model in our setting we write
$$ f(z)=\prod_{j\ge 1} \frac{1}{(1-z^j)^{b_j}},\ \text{where }\ b_j=\left\{
                                                                      \begin{array}{ll}
                                                                        1, & \hbox{if}\ j=b(mod\ a) \\
                                                                        0, & \hbox{otherwize}
                                                                      \end{array}
                                                                    \right.
$$
Thus, the sequence $b_j,\ j\ge 1$ of  weights has gaps. However,
$$  \sum_{1\le j\le n,q\!\!\not| j} b_j= O(n),\ \text{for\ all\ integers   }\ q>1,$$
by virtue of the fact that $gcd(a,b)=1$. Hence, \refm[njkv] holds, which secures the normal local limit theorem. The other ingredients of the asymptotic formula for the model are as follows.
$$D(s)=\zeta(s+1)D_b(s),\ D_b(s)= \sum_{k\ge 0}\frac{1}{(ak+b)^s}=a^{-s}\zeta(s,b/a),$$
where $\zeta(s,b/a),\ s\neq 1$ is the Hurwitz zeta- function, which has one simple pole $\rho_1=1$ with residue $A_1=1$;
$$h_1= a^{-1}\zeta(2)=a^{-1}\frac{\pi^2}{6},\ A_0= \zeta(0,b/a)
=1/2- b/a ,$$$$Q= h_0=-(\log a)(1/2- b/a)+ \log \Gamma(b/a) -1/2 \log(2\pi),$$$$ \kappa=1/2(-3/2 +1/2-b/a)=-1/2(1+ b/a)=-\frac{a+b}{2a}. $$
 The above  data lead to the following asymptotic formula:
$$c_n\sim \Gamma(\frac{b}{a})(\pi^{b/a-1})
(2^{1+ b/2a})(3^{-b/2a})n^{-\frac{a+b}{2a}}\exp\Big(\pi\sqrt{\frac{2n}{3a}}\Big).$$

This  formula was  obtained by Ingham in 1941(\cite{ing}), as a consequence of a Tauberian theorem for the growth of the generating function $f(z),$ as $z\to 1^-$.
Not long ago Kane (\cite{kan} derived two sided bounds on $c_n$ for the model considered, using elementary analytic tools.

\end{document}